# IMPLICATIVE FILTERS IN QUASI ORDERED *RL*-WAJSBERG ALGEBRAS

**Abstract** This article introduces the idea of implicative filters in quasi ordered RL-Wajsberg algebras and uses examples to explore some of its features.

## 1 INTRODUCTION

S. Bonzia [1] presented the idea of residual relational systems organized under a quasi-order relation in 2018. Ward and Dilworth introduced residuated lattices [6]. The idea of Wajsberg algebras was first suggested by Mordchaj Wajsberg in [7]. The terms filter and implicative filter were first used in lattice implication algebra by Y. Xu and K. Qin [13].

In this study, we introduce the concepts of implicative filters in Wajsberg algebra with a quasi-ordered residuated lattice (*RL)*. We also discover some of their associated qualities. We examine a few comparable conditions with examples before wrapping up.

## 2 MAIN RESULTS

### 2.1 IMPLICATIVE FILTER

In this part, we introduce the idea of implicative filters in Wajsberg algebra with a quasi-ordered residuated lattice. We also look into several implicative filter characteristics using examples.

**Definition 2.1.1.** Consider a Wajsberg algebra on a quasi-ordered residue lattice. If a subset $\mathcal{M}$ of $\wp$ meets the following axioms, it is said to be an implicative filter of $\wp$ for all $r, p, k \in \wp$.

(i) $1 \in \mathcal{M}$

(ii) If $((r \in \mathcal{M} \wedge r \leq p) \Longrightarrow p \in \mathcal{M})$

(iii) $r \twoheadrightarrow (p \twoheadrightarrow k) \in \mathcal{M} \wedge r \twoheadrightarrow p \in \mathcal{M} \Longrightarrow r \twoheadrightarrow k \in \mathcal{M}$.

**Proposition 2.1.2** In a quasi-ordered *RL* - Wajsberg algebra, every implicative filter is a filter.

***Proof.*** In a quasi-ordered *RL* - Wajsberg algebra, let $\mathcal{M}$ be the implicative filter, and let $r \twoheadrightarrow p \in \mathcal{M}$ and $r \in \mathcal{M}$ for all$r, p, k \in \wp$.

Replace $k$ by $p$ in definition 2.1.1.

Then, we have $r \twoheadrightarrow (p \twoheadrightarrow k) = r \twoheadrightarrow (p \twoheadrightarrow p)$

$$= r \twoheadrightarrow 1$$

$$= r \in \mathcal{M}.$$

Also, $1 \twoheadrightarrow (r \twoheadrightarrow p) \in \mathcal{M}$ and $1 \twoheadrightarrow r \in \mathcal{M}$. Then, $1 \twoheadrightarrow p \in \mathcal{M}$.

That implies, $p \in \mathcal{M}$.

Hence, in a quasi-ordered *RL* - Wajsberg algebra, every implicative filter is a filter.

**Proposition 2.1.3.** Let $\mathcal{M}$ be implicative filter in a quasi-ordered *RL* - Wajsberg algebra such that $r \twoheadrightarrow \big(p \twoheadrightarrow (p \twoheadrightarrow k)\big) \in \mathcal{M}$ and $r \in \mathcal{M}$ imply $r \twoheadrightarrow k$ for all $r, p, k \in \wp$. Then $\mathcal{M}$ is an implicative filter of $\wp$.

***Proof.*** Let $r \twoheadrightarrow \big(p \twoheadrightarrow (p \twoheadrightarrow k)\big) \in \mathcal{M}$ and $r \twoheadrightarrow p \in \mathcal{M}$ for all $r, p, k \in \wp$.

$r \twoheadrightarrow (p \twoheadrightarrow k) = p \twoheadrightarrow (x \twoheadrightarrow k) \leq (r \twoheadrightarrow p) \twoheadrightarrow \big(r \twoheadrightarrow (r \twoheadrightarrow k)\big)$

$(r \twoheadrightarrow p) \twoheadrightarrow \big(r \twoheadrightarrow (r \twoheadrightarrow k)\big) \in \mathcal{M}$

Since, $r \twoheadrightarrow p$, we have $r \twoheadrightarrow k \in \mathcal{M}$.

Hence, $\mathcal{M}$ is aimplicative filter of $\wp$.

**Proposition 2.1.4.** Let $\wp$ be a quasi-ordered *RL* - Wajsberg algebra, $\mathcal{M} \subseteq \wp$.Then, the following statements holds.

$\mathcal{M}$ is an implicative filter and for any $r, p \in \wp, r \twoheadrightarrow (r \twoheadrightarrow p) \in \mathcal{M}$ implies $r \twoheadrightarrow p \in \mathcal{M}$

***Proof.***

For any $r, p \in \wp$,if $r \twoheadrightarrow (r \twoheadrightarrow p) \in \mathcal{M}$, since $r \twoheadrightarrow r = 1 \in \mathcal{M}$, from (ii) of definition 2.1.1 we have $r \twoheadrightarrow p \in \mathcal{M}$.

**Proposition 2.1.5.** Let $\wp$ be a quasi-ordered *RL* - Wajsberg algebra, $\mathcal{M}_1$ and $\mathcal{M}_2$ are any two implicative filters of $\wp, \mathcal{M}_1 \subseteq \mathcal{M}_2$. If $\mathcal{M}_1$ is a implicative filter, so is $\mathcal{M}_2$.

***Proof.***

Suppose that, $r \twoheadrightarrow (r \twoheadrightarrow p) \in \mathcal{M}_2$, we only to prove $r \twoheadrightarrow p \in \mathcal{M}_2$.

Now,$r \twoheadrightarrow \Big(r \twoheadrightarrow \big(\big(r \twoheadrightarrow (r \twoheadrightarrow p)\big) \twoheadrightarrow p\big)\Big) = \big(r \twoheadrightarrow (r \twoheadrightarrow p)\big) \twoheadrightarrow \big(r \twoheadrightarrow (r \twoheadrightarrow p)\big) = 1 \in \mathcal{M}_1$

It follows that $r \twoheadrightarrow ((r \twoheadrightarrow (r \twoheadrightarrow p)) \twoheadrightarrow p \in \mathcal{M}_1 \subseteq \mathcal{M}_2$.

That is, $(r \twoheadrightarrow (r \twoheadrightarrow p)) \twoheadrightarrow (r \twoheadrightarrow p) \in \mathcal{M}_2$ and hence $r \twoheadrightarrow p \in \mathcal{M}_2$.

Hence $\mathcal{M}_2$ is a implicative filter of $\wp$.

**Definition 2.1.6.** Let $\wp$ be a quasi-ordered *RL* - Wajsberg algebra, $\psi \in \wp$. The interval $[\psi, 1]$ defined as $[\psi, 1] = \{r/r \in \wp, \psi \leq r\}$ of $\wp$ denoted as $\Upsilon(\psi)$.

**Proposition 2.1.7.** Let $\wp$ be a quasi-ordered RL- Wajsberg algebra, $\psi \in \wp$, then $\{1\}$ is an implicative filter of $\wp$ if and only if $\Upsilon(\psi)$ is an implicative filter of $\wp$ for any $\psi \in \wp$.

***Proof.*** Suppose that $\{1\}$ is an implicative filter of $\wp$. For any $\psi \in$ ,$1 \in \wp$ is trivial. If $r \in \Upsilon(\psi)$ and $r \twoheadrightarrow p \in \Upsilon(\psi)$, then $\psi \leq r, \psi \leq r \twoheadrightarrow p$, that is $\psi \twoheadrightarrow r = 1 \in \{1\}$ and $\psi \twoheadrightarrow (r \twoheadrightarrow p) = 1 \in \{1\}$. It follows that $\psi \twoheadrightarrow p \in \{1\}, a \leq p$, and hence $p \in \Upsilon(\psi)$. Thus, $\Upsilon(\psi)$ is an implicative filter of .

Conversely, assume that $\Upsilon(\psi)$ is an implicative filter of for any $\psi \in \wp$. For any $r, p, k \in \wp$. If $r \twoheadrightarrow (p \twoheadrightarrow k) \in \{1\}$ and $r \twoheadrightarrow p \in \{1\}$, then $r \leq p \twoheadrightarrow k, r \leq p$, it follows that $r \leq k$ because $\Upsilon(r)$ is an implicative filter, hence $r \twoheadrightarrow k = 1 \in \{1\}$. Hence, $\{1\}$ is an implicative filter of $\wp$.

**Proposition 2.1.8.** Let $\wp$ be a quasi-ordered *RL*- Wajsberg algebra and $\mathcal{M}$ be a subset of $\wp$ satisfies the condition $(r \in \mathcal{M} \wedge r.p) \Longrightarrow p \in \mathcal{M}$ for all $r, p, k \in \wp$.

***Proof.***

We know that $\psi \twoheadrightarrow r.r$ and $r.\psi \twoheadrightarrow r$ for all $r \in \wp$

From definition (ii) of 2.1.1 we have, $(r \in \mathcal{M} \wedge r.p) \Longrightarrow p \in \mathcal{M}$ for all $r, p, k \in \wp$.

**Proposition 2.1.9.** Let $\wp$ be a quasi-ordered *RL* - Wajsberg algebra, $\widehat{\wp} \subseteq \wp$.The following statements are equivalent for all $\tau, \varsigma, k \in \wp$.

(i) $\widehat{\wp}$ is implicative filter

(ii) $\widehat{\wp}$ is an implicative filter and for any $\tau, k \in \wp, \tau \twoheadrightarrow (\tau \twoheadrightarrow k) \in \widehat{\wp}$ implies $\tau \twoheadrightarrow k \in \widehat{\wp}$

(iii) $\widehat{\wp}$ is an implicative filter and for any $\tau, \varsigma, k \in \wp, \tau \twoheadrightarrow (k \rightarrow k) \in \widehat{\wp}$ implies $(\tau \twoheadrightarrow p) \twoheadrightarrow (\tau \twoheadrightarrow k) \in \widehat{\wp}$

(iv) $1 \in \widehat{\wp}$ and for any $\tau, \varsigma, k \in \wp, k \twoheadrightarrow (\tau \twoheadrightarrow (\tau \twoheadrightarrow k)) \in \widehat{\wp}$ and $k \in \widehat{\wp}$ imply $\tau \twoheadrightarrow k \in \widehat{\wp}$.

***Proof.***

**(i) ⇒ (ii)**

For any $\tau, k \in \wp$,if $\tau \twoheadrightarrow (\tau \twoheadrightarrow k) \in \widehat{\wp}$, since $\tau \twoheadrightarrow \tau = 1 \in \widehat{\wp}$, from (ii) of definition 2.1.1 we have $\tau \rightarrow k \in \widehat{\wp}$.

**(ii) ⇒ (iii)**

Assume that (ii) holds. For any $\tau, \varsigma, k \in \wp$, suppose $\tau \twoheadrightarrow (k \twoheadrightarrow k) \in \widehat{\wp}$ , implies

$\tau \twoheadrightarrow (\varsigma \twoheadrightarrow k) \leq \tau \twoheadrightarrow \big((\tau \twoheadrightarrow \varsigma) \twoheadrightarrow (\tau \twoheadrightarrow k)\big)$.

we get $\tau \twoheadrightarrow \Big(\tau \twoheadrightarrow \big((\tau \twoheadrightarrow \varsigma) \twoheadrightarrow k\big)\Big) = \tau \twoheadrightarrow \big((\tau \twoheadrightarrow \varsigma) \twoheadrightarrow (\tau \twoheadrightarrow k)\big) \in \widehat{\wp}$.

$\tau \twoheadrightarrow \big((\tau \twoheadrightarrow \varsigma) \twoheadrightarrow k\big) = (\tau \twoheadrightarrow \varsigma) \twoheadrightarrow (\tau \twoheadrightarrow k) \in \widehat{\wp}$

**(iii) ⇒(iv)**

Assume that (iii) holds and we prove (iv). Since, $1 \in \widehat{\wp}$ is trivial. For any $r, p, k \in \wp$, suppose $\varsigma \twoheadrightarrow (\tau \twoheadrightarrow (\tau \twoheadrightarrow k)) \in \widehat{\wp}$ and $k \in \widehat{\wp}$,then $\tau \twoheadrightarrow (k \twoheadrightarrow \varsigma) \in \widehat{\wp}$ .

Hence, we have $\tau \twoheadrightarrow \varsigma = 1 \twoheadrightarrow (\tau \twoheadrightarrow \varsigma) = (\tau \twoheadrightarrow \tau) \twoheadrightarrow (\tau \twoheadrightarrow \varsigma) \in \widehat{\wp}$.

**(iv) ⇒(i)**

Suppose $\tau \in \widehat{\wp}$ and $\tau \twoheadrightarrow \varsigma \in \widehat{\wp}$ then we get $\tau \twoheadrightarrow \big(1 \twoheadrightarrow (1 \twoheadrightarrow \varsigma)\big) = \tau \twoheadrightarrow \varsigma \in \widehat{\wp}$, it follows that $\varsigma = 1 \twoheadrightarrow \varsigma \in \widehat{\wp}$ and hence $\widehat{\wp}$ is an implicative filter.

For any $\tau, \varsigma, k \in \wp, \tau \twoheadrightarrow (\varsigma \twoheadrightarrow k) \in \widehat{\wp}$ and $\tau \twoheadrightarrow \varsigma \in \widehat{\wp}$,

Now, $\left(r \twoheadrightarrow (\varsigma \twoheadrightarrow k)\right) \twoheadrightarrow \left((\tau \twoheadrightarrow \varsigma) \twoheadrightarrow \left(\tau \twoheadrightarrow (\tau \twoheadrightarrow \varsigma)\right)\right)$

$$= (\tau \twoheadrightarrow \varsigma) \twoheadrightarrow ((\tau \twoheadrightarrow (\varsigma \twoheadrightarrow k)) \twoheadrightarrow (\tau \twoheadrightarrow (\tau \twoheadrightarrow k)))$$

$$= (\tau \twoheadrightarrow \varsigma) \twoheadrightarrow ((\tau \twoheadrightarrow (\tau \twoheadrightarrow \varsigma)) \twoheadrightarrow (\tau \twoheadrightarrow (\tau \twoheadrightarrow \varsigma)))$$

$$= (\tau \twoheadrightarrow \varsigma) \twoheadrightarrow (\tau \twoheadrightarrow (\varsigma \vee (\tau \twoheadrightarrow k)))$$

$$= (\tau \twoheadrightarrow \varsigma) \twoheadrightarrow \left((\tau \twoheadrightarrow k) \vee \left(\tau \twoheadrightarrow (\tau \twoheadrightarrow \varsigma)\right)\right) = 1 \in \hat{\wp}.$$

Hence $(\tau \twoheadrightarrow \varsigma) \twoheadrightarrow (\tau \twoheadrightarrow (\tau \twoheadrightarrow \varsigma)) \in \hat{\wp}$and hence $\tau \twoheadrightarrow k \in \hat{\wp}$ by $\tau \twoheadrightarrow \varsigma \in \hat{\wp}$ and (iv).

**Proposition 2.1.10.** If a non-empty subset $\mathcal{M}$ of a quasi-ordered *RL* - Wajsberg algebra $\wp$, satisfies conditions for all $\tau, \varsigma, k \in \wp$

(i) $(\tau \in \mathcal{M} \wedge \tau \odot \varsigma) \Rightarrow \varsigma \in \mathcal{M}$
(ii) $(\tau \in \mathcal{M} \wedge \tau \twoheadrightarrow \varsigma \in \mathcal{M}) \Rightarrow \varsigma \in \mathcal{M}$
(iii) $(\tau \twoheadrightarrow (\tau \twoheadrightarrow \varsigma) \in \mathcal{M}) \Rightarrow \tau \in \varsigma \in \mathcal{M}$

Then $\mathcal{M}$ is an implicative filter in $\wp$.

***Proof.***

***To prove***: $\mathcal{M}$ is an implicative filter in $\wp$

It is enough to show that $(\tau \twoheadrightarrow \varsigma \twoheadrightarrow (\varsigma \twoheadrightarrow k)) \in \mathcal{M}) \wedge \ \tau \in \mathcal{M} \Rightarrow \varsigma \twoheadrightarrow k \in \mathcal{M}$

Let $r, \varsigma, k \in \wp$ such that $(\tau \twoheadrightarrow \varsigma \twoheadrightarrow (\varsigma \twoheadrightarrow k)) \in \mathcal{M})$ and $\tau \in \mathcal{M}$

Then, $(\varsigma \twoheadrightarrow (\varsigma \twoheadrightarrow k)) \in \mathcal{M})$ and$\tau \in \mathcal{M}$.

Hence, $\varsigma \twoheadrightarrow k \in \mathcal{M}$.

## 3. CONCLUSION

The principles of quasi-ordered RL - Wajsberg algebra have been put forth. We also covered a similar requirement that every implicative filter is an implicative filter. Additionally, we obtained some implicative filter aspects in quasi-ordered RL-Wajsberg algebra.